%%%%%%%%%%%%%%%%%%%%%%%%%%%%%%
% Version 1/23/96            %
%%%%%%%%%%%%%%%%%%%%%%%%%%%%%%
\magnification 1200
\input amstex
\documentstyle{amsppt}
%\NoRunningHeads
\NoBlackBoxes
%\hoffset 5.5truemm
%\vsize 7.5in
\topmatter
\title CR automorphisms of real analytic 
manifolds in complex space\endtitle
\rightheadtext{CR automorphisms of real analytic 
manifolds}
\leftheadtext{M.~S.~Baouendi, P.~Ebenfelt, and
L.~P.~Rothschild}
\author M. S. Baouendi\footnote{{Partially supported by National
Science Foundation Grant DMS 
95-01516.\hfill\break}},  P.
Ebenfelt\footnote{Supported by a grant from the Swedish
Natural Science Research Council.\newline}, and Linda Preiss
Rothschild$^1$
\endauthor

%\abstract \endabstract
%\keywords \endkeywords   
%\subjclass \endsubjclass
\address Department of Mathematics, University of
California at San Diego, La Jolla, CA 92093\endaddress
\email sbaouendi\@ucsd.edu, 
pebenfel\@euclid.ucsd.edu, lrothschild\@ucsd.edu\endemail
%\date{\number\year-\number\month-\number\day}\enddate
\loadeufm
%ndefine \p{\eufm p}
%\define \q{\eufm q}
%\define \m{\eufm m}
%\define \ass{\text{\rm Ass }}
%\define \sol{\text{\rm Sol}}

\define \im{\text{\rm Im }}
\define \re{\text{\rm Re }}

\define \bR{\Bbb R}
%\define \bR{{\bold R}}
\define \bC{\Bbb C}
%\define \bC{{\bold C}}
\define \bL{\Bbb L}
%\define \bL{{\bold L}}
\define \scrM{\Cal M}

\define \scrO{\Cal O}
\define \scrV{\Cal V}

%\define \ssneq{\subsetneqq}
%\define \ssneq{\subset}
\define \hol{\text{\rm hol}}
\define \aut{\text{\rm aut}}

\def\cnn {\bC^N}

\def \po {p_0}
\def\r {\rho}

\def\z {\zeta}

\def\dim {\text {\rm dim}}

\endtopmatter
\document 
\heading 0. Introduction \endheading

 In this paper we shall give sufficient conditions for
local CR diffeomorphisms between two real analytic
submanifolds of
$\bC^N$ to be determined by finitely many derivatives at
finitely many points. These conditions will also be shown
to be necessary in model cases.   We shall also show that
under the same conditions, the Lie algebra of the
infinitesimal CR automorphisms at a point is finite
dimensional.  

Let $M$ be a real analytic submanifold of $\bC^N$.  For
$p\in M$ a {\it CR vector} at
$p$ is a vector of the form $\sum_{j=1}^N c_j {\partial
\over
\partial {\overline Z}_j}$, $c_j \in \bC$, tangent to $M$ at
$p$. If
$M'$ is another submanifold of $\bC^N$,  a mapping
$F:M
\to M'$
 is called {\it CR} if for any $p \in M$ the pushforward
$F_* X$ of any CR vector $X$ on $M$ at $p$ is a CR vector
of $M'$ at $F(p)$.   In particular, the restriction to
$M$ of a germ of a holomorphic diffeomorphism $H$ from
$\bC^N$ to itself is a CR map from $M$ to its image.

As in [BER1] (see Stanton [St1] for the case of a
hypersurface), we shall say that a real submanifold of
$\bC^N$ is {\it holomorphically nondegenerate} if there is
no germ of a nontrivial vector field $\sum_{j=1}^N c_j(Z)
{\partial \over
\partial Z_j}$, with $c_j(Z)$ holomorphic, tangent to
$M$.  If $M$ is 
holomorphically nondegenerate there is an integer
$l(M)$, with $0 \le l(M) < N$, called the  Levi
number of
$M$
 (see
\S 1) which  measures the holomorphic nondegeneracy of
$M$. If
$M$ is a Levi-nondegenerate hypersurface then $l(M) =
1$.   A connected real analytic submanifold is {\it
minimal almost everywhere} if there is no germ of a
holomorphic function whose restriction to $M$ is a
nonconstant real-valued function. This coincides with the
notion of being minimal at most points in the sense of
Tumanov [Tu1].  If
$M$ is a hypersurface which is holomorphically nondegenerate,
then
$M$ is minimal almost everywhere.  

The following uniqueness result is one of the main
theorems of this paper.      

\proclaim {Theorem 1} Let $M\subset\bC^N$  be a
connected, real analytic, holomorphically nondegenerate
submanifold of codimension
$d$ and Levi number $l(M)$ such that $M$ is minimal almost
everywhere. Then for all $p \in M$ outside a proper real
analytic subvariety of $M$ the following holds.
  If $M^\prime\subset \bC^N$ is another real analytic
submanifold with 
$\dim_\bR M^\prime=\dim_\bR M$, and
$F,G$ are smooth germs at $p$ of CR diffeomorphisms  of
$M$ into
$M^\prime$ such that in some local  coordinates $x$ on $M$
$$
\frac{\partial^{|\alpha|}F}{\partial x^\alpha}(p)=
\frac{\partial^{|\alpha|}G}{\partial x^\alpha}(p)
\tag0.1
$$ for all $|\alpha|\leq (d+1)l(M)$,  then $F\equiv G$. 
\endproclaim

\proclaim {Corollary} Let $M$ be as in Theorem 1.  Then for all $p \in M$ outside a proper real
analytic subvariety of $M$ the following holds. If $H$ is a
germ at $p$ of a local biholomorphism of
$\bC^N$ mapping $M$ into itself and fixing $p$,  with
$$
\frac{\partial H_j}{\partial Z_k}(p)=\delta_{jk},\ \ \  \
\frac{\partial^{|\alpha|}H_j}{\partial Z^\alpha}(p) =0, \
\ 1\le j,k \le N, \  2 \le |\alpha| \leq (d+1)l(M), 
\tag0.2 
$$ then $H$ is the identity map on $M$. 
\endproclaim  In fact, Theorem 1 also follows from the
statement of the Corollary.
 In  case $M$ is a Levi-nondegenerate hypersurface (i.e.
$d=1$ and $l(M)=1$) Theorem 1 reduces to the
result of Chern-Moser [CM] that a germ of a CR
diffeomorphism is uniquely determined by its derivatives
of order $\le 2$ at a point.  Generalizations of this
result for Levi nondegenerate manifolds of higher
codimension were later given by Tumanov-Henkin [TH],
Tumanov [Tu2]. More precise results for
Levi nondegenerate hypersurfaces have been given by
Beloshapka [Be] and Loboda [L].

A smooth real vector field $X$ defined in a neighborhood
of
$p$ in $M$  is an {\it infinitesimal holomorphism } if the
 local 1-parameter group of diffeomorphisms $\exp tX$ for $t$
small extend to a local 1-parameter group of
biholomorphisms of $\bC^N$.  More generally,
$X$ is called an {\it infinitesimal CR  automorphism } if 
the
$\exp tX$ are CR diffeomorphisms.  We denote by
$\hol(M,p)$ (resp. $\aut(M,p)$) the Lie algebra
generated by the infinitesimal holomorphisms (resp.
infinitesimal CR automorphisms).  Since every local
biholomorphism preserving $M$ restricts to a CR
diffeomorphism of $M$ into itself, it follows that
$\hol(M,p)
\subset\aut(M,p)$. It follows from the work of Tanaka
[Ta] that $\hol(M,p)$ is a finite dimensional vector space if
$M$ is a real analytic Levi nondegenerate hypersurface.
Recently Stanton [St2] proved that if
$M$ is a real analytic hypersurface, $\hol(M,p)$ is a finite
dimensional real vector space for any 
$p \in M$ if and only if $M$ is holomorphically
nondegenerate.  In this paper we prove  more general
results for any real analytic CR submanifold. Recall that a
real analytic submanifold $M$ of $\bC^N$ of codimension
$d$ is {\it CR} if  $M$ is locally defined by the
vanishing of
$d$ real valued real analytic functions $\rho_1, \ldots,
\rho_d$, with linearly independent differentials, such that
the linear span of the complex differentials
$\partial
\rho_1,
\ldots, \partial \rho_d $ is of constant dimension. 

\proclaim {Theorem 2}
 Let $M\subset\bC^N$ be a real analytic, connected CR
submanifold. If $M$ is holomorphically nondegenerate, and
minimal almost everywhere then 
$$
\dim_\bR \aut(M,p)<\infty\tag0.3
$$ for all $p \in M$.  
\endproclaim

 Theorems 1 and 2 are optimal in the sense that
holomorphic nondegeneracy is necessary for the
conclusions of Theorems 1 and 2 and that the condition
that $M$ is minimal almost everywhere is necessary in
model cases.  We have the following result.  

\proclaim {Theorem 3} Let $M\subset\bC^N$ be a connected
real analytic CR submanifold. 
\roster
\item "(i)" If $M$ is holomorphically
degenerate, then for any $p \in M$ and any integer $K >
0$ there exist local biholomorphisms $F$ and $G$ near $p$
mapping $M$ into itself and fixing $p$ such that
 $$
\frac{\partial^{|\alpha|}F}{\partial Z^\alpha}(p)=
\frac{\partial^{|\alpha|}G}{\partial Z^\alpha}(p)
\tag0.4
$$ for all $|\alpha|\leq K$,  but $F \not\equiv G$ on
$M$. Furthermore, {\rm dim}$_\bR\hol(M,p) =
\infty$.
\item "(ii)"If $M$ is holomorphically nondegenerate but
nowhere minimal then for $p$ in an open dense set
in $M$ either {\rm dim}$_\bR\hol(M,p) =
\infty$ or {\rm dim}$_\bR\hol(M,p) =0$. 
\item "(iii)" If
$M$ is defined by the vanishing of weighted homogeneous
polynomials, and nowhere minimal then for any $p \in M$ and
any integer
$K > 0$ there exist local biholomorphisms $F$ and $G$ near
$p$ mapping $M$ into itself and fixing $p$ such that (0.4)
holds for all $|\alpha|\leq K$,  but $F \not\equiv G$ on
$M$. Also
{\rm dim}$_\bR\hol(M,p) = \infty$.
\endroster   
\endproclaim

The paper is organized as follows.  In \S 1 we recall some
results about holomorphically nongenerate manifolds and Segre
sets of a generic manifold.  In \S 2 we prove a result on
uniqueness of CR diffeomorphisms, of which Theorem 1 is a
consequence.  In \S 3 we study infinitesimal CR
automorphisms and prove Theorem 2.  In \S 4 we prove Theorem
3 and also construct an example to show that one can have
{\rm dim}$_\bR \hol(M,p) = 0$ even if $M$ is nowhere
minimal. In \S 5 we make some remarks concerning the group
$G_p$ of local biholomorphisms leaving $M$ invariant and
fixing a point $p \in M$.

\heading 1. Preliminaries\endheading
\subhead 1.1. Holomorphic nondegeneracy and
$k$-nondegeneracy of real analytic CR
submanifolds\endsubhead Let $M\subset\bC^N$ be a
connected real analytic CR submanifold.  We denote by 
$\scrV$ the largest holomorphic submanifold in $\bC^N$
containing $M$ with minimum dimension. As in [BER1] we
call 
$\scrV$ the {\it intrinsic complexification} of $M$. If
 $\scrV$ is all of $\bC^N$, then
$M$ is called {\it generic}. Thus any CR submanifold is
generic when considered as a submanifold of its intrinsic
complexification.  Recall that
$M$ is called {\it holomorphically degenerate} at $p_0
\in M$ if there is a germ of a  holomorphic vector field
at
$p_0$ which is tangent to $M$ but not trivial  (i.e. not
identically 0) on $M$. 
 It can be easily checked that  $M$ is holomorphically
degenerate at
$p_0$ as a submanifold of $\bC^N$ if and only if it is
holomorphically degenerate at $p_0$ as a submanifold of
$\scrV$.

We say that $M$ is {\it holomorphically nondegenerate} 
if it is not holomorphically degenerate at any point.  In
fact it is proved in [BER1] that if $M$ is holomorphically
degenerate at one point, it is holomorphically degenerate
at every point.

Suppose $M$ is a real analytic generic submanifold of
$\bC^N$ and that
$M$ is defined near $p_0\in M$ by $\rho(Z,\bar Z)=0$,
where $\rho=(\rho_1,...,
\rho_d)$ are real valued real analytic functions with
$\rho(p_0,\bar p_0)=0$ and
$\partial\rho_1\wedge...\wedge\partial\rho_d\neq0$ near
$p_0$. Let $L=(L_1,...,L_n)$ be a basis for the CR vector
fields on
$M$ near $p_0$. For any  multi-index
$\alpha$ put  $L^\alpha = L_1^{\alpha_1}\ldots 
L_n^{\alpha_n}$. Introduce, for
$j=1,...,d$ and any  multi-index
$\alpha$, the vectors
$$ V_{j\alpha}(Z,\bar Z)=L^\alpha\rho_{jZ}(Z,\bar
Z),\tag1.1.1$$
 where $\rho_{jZ}$ denotes the gradient of
$\rho_j$ with respect to $Z$.  We say that the generic
real analytic submanifold
$M$ is {\it $k$-nondegenerate} at $p_0$ if $k$ is the
smallest positive integer for which the span of the
vectors
$V_{j\alpha} (p_0,\bar p_0)$, for
$j=1,...,d$ and $|\alpha|\leq k$, equals $\bC^N$. This
definition is independent of the coordinate system used,
the defining equations of $M$, and the choice of basis
$L$.  We say that a real analytic CR submanifold $M$ is
$k$-nondegenerate at $p_0\in M$ if $M$ is
$k$-nondegenerate at $p_0$ as a generic submanifold  of
its intrinsic complexification $\scrV$.

If $M$ is a connected CR submanifold of $\bC^N$ and $\scrV$
its intrinsic complexification, then the {\it CR dimension}
of $M$ is the nonnegative integer defined by
$$ \dim_\bR M = \dim_\bC \scrV + CR\ \dim M.
$$ 
 The following proposition is  in [BER1,
Proposition 1.3.1] for generic manifolds.  Its extension to
CR manifolds is immediate. 

\proclaim{Proposition 1.1.1} Let $M\subset\bC^N$ be a
connected real analytic CR submanifold of CR dimension
$n$. Then  the following are equivalent.
\roster
\item"(i)" $M$ is holomorphically nondegenerate.
\item"(ii)" There exists $p_1\in M$ and $k>0$ such that
$M$ is
$k$-nondegenerate at $p_1$. 
\item"(iii)" There exists $V$, a proper real analytic
subset of $M$, and an integer $l=l(M)$, $0\leq l(M)\leq
n$, such that $M$ is $l$-nondegenerate at every $p\in
M\setminus V$.
\endroster\endproclaim

The number $l(M)$ given in (iii) is called the Levi
number of $M$.

\subhead 1.2. The Segre sets\endsubhead In this section,
we introduce the Segre sets of a generic real analytic
submanifold in $\bC^N$ and recall some of their
properties. We refer the reader to the paper [BER1] for a
more detailed account (including proofs of the main
results) of these sets. Let
$M$ denote a generic  real analytic submanifold in some
neighborhood
$U\subset \bC^N$ of $p_0\in M$. Let $\r = (\r_1,\ldots
\r_d) $ be defining functions as above, and choose
holomorphic coordinates $Z=(Z_1,\ldots,Z_N)$ vanishing at
$\po$.  Embed 
$\bC^N$ in
$\bC^{2N}=\bC^N_Z\times\bC^N_\zeta$ as the real plane
$\{(Z,\zeta)\in\bC^{2N}\:\zeta=\bar Z\}$. Let us denote by
$\text{pr}_Z$ and $\text{pr}_\zeta$ the projections of
$\bC^{2N}$ onto 
$\bC^N_Z$ and $\bC^N_\zeta$, respectively. The natural
anti-holomorphic involution
$\sharp$ in $\bC^{2N}$ defined by
$$ ^\sharp(Z,\zeta)=(\bar\zeta,\bar Z)\tag1.2.1
$$ 
 leaves the plane
$\{(Z,\zeta)\:\zeta=\bar Z\}$ invariant. This  involution
induces the usual anti-holomorphic involution in $\bC^N$
by
$$
\bC^N\ni
Z\mapsto\text{pr}_\zeta(^\sharp\text{pr}_Z^{-1}(Z))=\bar
Z\in\bC^N.
\tag1.2.2
$$ Given a set $S$ in $\bC^N_Z$ we denote by
$^* S$ the set in
$\bC^N_\zeta$ defined by
$$
^*S=\text{pr}_\zeta(^\sharp\text{pr}_Z^{-1}(S))=\{\zeta\:\bar\zeta\in
S\}.
\tag1.2.3
$$ By a slight abuse of notation, we use the same
notation for the corresponding  transformation taking
sets in $\bC^N_\zeta$ to sets in $\bC^N_Z$. Note that if
$X$ is a complex analytic set defined  near
$Z^0$ in some domain $\Omega
\subset\bC_Z^{N}$  by
$h_1(Z)=...=h_k(Z)=0$,  then $^* X$ is the complex  
analytic set in
$^*\Omega\subset\bC^N_\zeta$ defined near
$\zeta^0=\bar Z^0$ by
$\bar h_1(\zeta)=...=\bar h_k(\zeta)=0$. Here, given a
holomorphic function $h(Z)$ we use the notation 
$
\bar h(Z)=\overline{h(\bar Z)}$. 

Denote by
$\scrM\subset\bC^{2N}$ the complexification of $M$ given
by
$$
\scrM=
\{(Z,\zeta)
\in\bC^{2N}
\:\rho(Z,\zeta)=0\}.\tag1.2.4$$  This is a complex
submanifold of codimension $d$ in some neighborhood of 
$0$ in $\bC^{2N}$. We choose our neighborhood
$U$ in $\bC^N$ so small that
$U\times{}^*U\subset\bC^{2N}$ is contained in the
neighborhood where
$\scrM$ is a manifold. Note that
$\scrM$ is invariant under the involution 
$\sharp$ defined in (1.2.1). 

We associate to
$M$ at $p_0$ a sequence  of germs of sets
$N_0,N_1,...,N_{j_0}$ at
$p_0$ in $\bC^N$---{\it the Segre sets} of $M$ at
$p_0$---defined as follows.  Put $N_0=\{p_0\}$ and
define the consecutive sets inductively (the number
$j_0$ will be defined later) by
$$ N_{j+1}=\text{pr}_Z\left(\scrM\cap
\text{pr}_\zeta^{-1}\left(^*N_j\right)
\right)=\text{pr}_Z\left(\scrM\cap{}^\sharp\text{pr}_Z^{-1}\left(N_j\right)
\right).\tag1.2.5
$$  Here, and in what follows, we  identify a  germ
$N_j$ with some representative of it. These sets are, by
definition, invariantly defined and they  arise naturally
in the study of mappings between submanifolds (see later
sections in this paper, and [BER1]).

The sets
$N_j$ can be described in terms of the defining equations
$\rho(Z,\bar Z)=0$ (see [BER1, \S 2.2]), e.g.
$$  N_1=\{Z\:\rho(Z,0)=0\}\tag1.2.6
$$  and
$$  N_2=\{Z\:\exists
\zeta^1\:\rho(Z,\zeta^1)=0\,,\,\rho(0,\zeta^1)
=0\}.\tag1.2.7$$ We have the inclusions
$$ N_0\subset N_1\subset...\subset
N_j\subset...\tag1.2.8$$   and  $j_0$ is the
largest number such that the first 
$j_0$ inclusions in \thetag{1.2.8} are strict. (The Segre
sets stabilize after that, and $N_{j_0+1}=N_{j_0+2}=...$
.) It is shown in [BER1] that, in suitable coordinates
$(z,w)\in\bC^n\times\bC^d=\bC^N$  (so-called normal
coordinates), the Segre set $N_j$, for $j=1,...,j_0$, can
also be defined as images of  certain holomorphic mappings
$$
\bC^n\times\bC^{(j-1)n}\ni(z,\Lambda)\mapsto
(z,v^j(z,\Lambda))\in\bC^N.\tag 1.2.9
$$ Thus, we can define the generic dimension $d_j$ of
$N_j$ as the  generic rank of the mapping \thetag{1.2.9}. 

So far we have only considered generic  submanifolds. If
$M$ is a real analytic CR submanifold of $\bC^N$, then
$M$ is generic as a submanifold of its intrinsic
complexification $\Cal V$ (see \S1.1). The Segre sets of
$M$ at a point
$p_0\in M$ can be defined  as subsets of
$\bC^N$  by the process described at the beginning of
this subsection (i.e. by
\thetag{1.2.5}) just as for generic submanifolds or they
can be defined as subsets of $\Cal V$ by identifying
$\Cal V$ near $p_0$ with $\bC^K$ and considering
$M$ as a generic submanifold of $\bC^K$. It can be shown
that these definitions are equivalent.

The main properties concerning the Segre sets  that we
shall use in this paper are summarized in the following
theorem.

\proclaim{Theorem 1.2.1} Let $M$ be a  real analytic CR
submanifold in
$\cnn$, and let $p_0\in M$. 
\roster
\item"(a)" Denote by $W$ the CR orbit of $p_0$ (i.e. the
Nagano leaf or, equivalently, the  CR submanifold
of $M$ of smallest dimension, with  the same CR dimension as
$M$, through
$p_0$) and by $X$ its intrinsic  complexification (see
\S1.1). Then the maximal Segre set $N_{j_0}$ of $M$ at
$p_0$ is contained in $X$ and $N_{j_0}$ contains an open
subset of $X$, i.e. $d_{j_0}=
\dim_{\bC} X$. 
\item"(b)" There are holomorphic immersions
$Z_0(t_0),Z_1(t_1),...,Z_{j_0} (t_{j_0})$ defined near the
origin,
$$
\bC^{d_j}\ni t_j\mapsto Z_j(t_j)\in\bC^N,\tag1.2.10$$ and
holomorphic maps $s_0(t_1),...,s_{j_0-1}(t_{j_0})$,
$$
\bC^{d_j}\ni t_j\mapsto
s_{j-1}(t_j)\in\bC^{d_{j-1}},\tag1.2.11
$$ such that $Z_j(t_j)$ has rank $d_j$ near the origin,
$Z_j(t_j)\in N_j$, and such that
$$
\left(Z_j(t_j),\bar Z_{j-1}(s_{j-1}(t_j))\right)\in\Cal
M,\tag1.2.12
$$ for $j=1,...,j_0$.\endroster
\endproclaim  
\demo{Proof} Part (a) is contained in [BER1, Theorem
1.2.2], and the mappings in part (b) are constructed in
the paragraph following [BER1, Assertion 3.3.2].
\qed\enddemo
\noindent{\bf Remark 1.2.2.} The holomorphic immersion
$Z_j(t_j)$,
$j=0,1,...,j_0$, in part (b) above provides a
parametrization of an open piece of $N_j$. However, this
piece of $N_j$ need not contain the point $p_0$. Indeed,
$N_j$ need not even be a manifold at $p_0$. \medskip

Recall that a CR submanifold $M$ is said to be minimal at
a point $p_0\in M$ if there is no proper CR submanifold
of $M$ through $p_0$ with the same CR dimension as $M$.
For a real analytic submanifold, this notion coincides
with the notion of finite type in the sense of
Bloom--Graham [BG]. One can check that if $M$ is connected
then $M$ is minimal almost everywhere, as defined in \S 0
if and only if $M$ is minimal at some point in the sense
described above. The following is an immediate consequence of
the theorem.
\proclaim{Corollary 1.2.3} Let $M$ be a  real analytic
generic submanifold in
$\cnn$ and $p_0\in M$. Then
$M$ is minimal at $p_0$  if and only if $d_{j_0} = N$ or,
equivalently, if and only if the maximal Segre set at $p_0$
contains an open subset of $\bC^N$.
\endproclaim

\heading 2. Uniqueness of CR diffeomorphisms\endheading
The main result here is the following, which implies
Theorem 1 as a special case.
\proclaim{Theorem 2.1} Let $M\subset\bC^N$ be a connected
real analytic, holomorphically nondegenerate CR
submanifold and let $d$ be the (real) codimension of $M$ in
its intrinsic complexificiation.  Suppose that there is a
point
$p\in M$ at which $M$ is minimal. For any $p_0\in M$ 
there exists a finite set of points
$p_1,...,p_k\in M$ such that if $M^\prime\subset \bC^N$
is another real analytic CR submanifold with 
$\dim_\bR M^\prime=\dim_\bR M$, and
$F,G$ are smooth  CR diffeomorphisms of $M$ into
$M^\prime$ such that 
$$
\frac{\partial^{|\alpha|}F}{\partial x^\alpha}(p_l)=
\frac{\partial^{|\alpha|}G}{\partial x^\alpha}(p_l)
\tag2.1
$$ for $l=1,...,k$, and $|\alpha|\leq (d+1)l(M)$,  then
$F\equiv G$ in a neighborhood of
$\po$ in $M$.  If
$M$ is minimal at
$p_0$, then one can take $k = 1$.  If, in addition, $M$ is
$l(M)$-nondegenerate at $\po$, then one may take $ p_1 =
\po$.
\endproclaim
\noindent{\bf Remarks:} \roster
\item"(i)" The condition \thetag{2.1} can be expressed by 
saying that the $(d+1)l(M)$-jets of the mappings coincide
at all the points
$p_1,...,p_k$. 
\item"(ii)" The choice of points $p_1,...,p_k$ can be
described as follows. Let $U_1,...,U_k$ be the components
of the set of minimal points of $M$ in 
$U$, a sufficiently small neighborhood of $p_0$ in $M$,
which have $p_0$ in their closure. For each $l=1,...,k$,
we may choose any $p_l$ from the dense open subset of
$U_l$ consisting of those points which are
$l(M)$-nondegenerate.
\endroster
\medskip

Before we prove Theorem 2.1 we need some preliminary
results. 

\proclaim{Proposition 2.2} Let $M,M^\prime\subset\bC^N$
be real analytic CR submanifolds, and $p_0\in M$. Assume
that $M$ is holomorphically nondegenerate and generic,
and that $M$ is
$l(M)$-nondegenerate at $p_0$. Let $H$ be a germ of a
biholomorphism of $\bC^N$ at $p_0$ such that $H(M)\subset
M^\prime$.  Then there are $\bC^N$ valued  functions
$\Psi^\gamma$, holomorphic in all of their arguments,
such that
$$
\frac{\partial^{|\gamma|} H}{\partial
Z^\gamma}(Z)=\Psi^\gamma\left(Z,\zeta,\bar H(\z), \ldots,
\frac{\partial^{|\alpha|}\bar H}{\partial
\zeta^\alpha}(\zeta),\ldots\right),
\tag2.2
$$ where $|\alpha|\leq l(M)+|\gamma|$, for all
multi-indices $\gamma$ and all points $(Z,\zeta)\in\scrM$
near $(p_0,\bar p_0)$. Moreover, the 
functions $\Psi^\gamma$  depend only on $M,M^\prime$ and 
$$
\frac{\partial^{|\beta|} H}{\partial Z^\beta}(p_0),\
\ |\beta|\leq l(M). \tag2.3
$$  \endproclaim
\demo{Proof} It suffices to prove \thetag{2.2} in any
coordinate system of the target space near
$p_0^\prime=H(p_0)$.  If we choose normal coordinates for
$M^\prime$ at $p_0^\prime$ then the proof is exactly the
same as the proof of Assertion 3.3.1 and subsequent
remarks in [BER1] (see also [BR1, Lemma 2.3]). \qed\enddemo

\proclaim{Proposition 2.3} Let $M\subset\bC^N$ be a real
analytic, holomorphically nondegenerate  CR submanifold and
$p_0\in M$. Let $X$ be the intrinsic complexification
of the CR orbit $W$ of $p_0$ in $M$ and $d$ the codimension
of
$W$ in $X$. Assume that
$M$ is
$l(M)$-nondegenerate at $p_0$.  If $H^1,H^2$ are germs of
biholomorphisms of $\bC^N$ at $p_0$ such that
$H^1(M), H^2(M)\subset M^\prime$, where $M^\prime$ is
another real analytic CR submanifold with $\dim_\bR
M^\prime=\dim_\bR M$, and such that
$$
\frac{\partial^{|\alpha|}H^1}{\partial Z^\alpha}(p_0)=
\frac{\partial^{|\alpha|}H^2}{\partial
Z^\alpha}(p_0),\quad\ |\alpha|
\leq (d+1)l(M),\tag2.4
$$ then $H^1|_X\equiv H^2|_X$.\endproclaim

\demo{Proof} The intrinsic complexification of the
CR orbit is contained in the intrinsic  complexification
of
$M$, and the notion of holomorphic nondegeneracy is
independent of the ambient space. Similarly, the notion of
$l(M)$-nondegeneracy is defined in the intrinsic
complexification.  Hence we  may reduce to the case where $M$
is generic;  we shall assume this for the rest of the proof. 
 
 Let $N_j$, $j=0,1,...,j_0$, be the Segre sets of
$M$ at
$p_0$, and let $Z_0(t_0),...,Z_{j_0}(t_{j_0})$ be the 
canonical parametrizations of the $N_j$'s and
$s_0(t_1),...,s_{j_0-1}$ $(t_ {j_0})$ the associated maps
so that
$$
\left(Z_{j+1}(t_{j+1}),\bar
Z_j(s_j(t_{j+1}))\right)\in\scrM,\tag2.5
$$ for all $j=0,...,j_0-1$ (see Theorem 1.2.1 (b)). In
view of 
\thetag{2.4} and Proposition  2.2,  there are functions
$\Psi^\gamma$ such that both $H^1$ and $H^2$ satisfy  the
identity \thetag{2.2} for $(Z,\zeta)\in\scrM$.
Substituting \thetag{2.5} with $j=0$  into this identity
and recalling that $Z_0(t_0)\equiv p_0$ (i.e. it is the 
constant map), we deduce that $H^1$ and $H^2$ as well as
all their derivatives are identical on the first Segre set
$N_1$. Note that since each $N_j$ is the holomorphic image
of a connected set, if two holomorphic functions agree on
an open piece, they agree on all of $N_j$.   By inductively
substituting
\thetag{2.5} into
\thetag{2.2} for $j=1,...,j_0-1$, we deduce that the
restrictions of the mappings $H^1$ and $H^2$, as well as all
their derivatives, to the maximal Segre set $N_{j_0}$ are 
identical.  The conclusion of the proposition now follows
from Theorem 1.2.1 (c), since $N_{j_0}$ contains an open
piece of $X$.
\qed\enddemo

We now proceed with the proof of Theorem 2.1.
\demo{Proof of Theorem 2.1} By a change of holomorphic
coordinates near $\po$ and by shrinking $M$, if
necessary,  we may assume that $M$ is a real analytic,
holomorphically nondegenerate, generic submanifold of
$\bC^K$, for some
$K\leq N$. 

Let $V$ be the set of points on $M$ at which
$M$ is not minimal. This is a real analytic subset and,
since $M$ is minimal at some point $p$, $V$ is also
proper. Denote by $U_1,...,U_k$ those components of
$M\setminus V$ that have $p_0$ in their closure. Clearly,
$k$ is a finite number since $M\setminus V$ is a
semi-analytic subset (a semi-analytic set is locally
finite in the sense that only a finite number of
components meet each compact set). Also, since $M$ is
holomorphically  nondegenerate, $M$ is
$l(M)$-nondegenerate outside a proper real analytic
subset. Pick $p_l\in U_l$, for $l=1,...,k$, such that $M$
is $l(M)$-nondegenerate at $p_l$. 

Since $M$ is minimal in $U_l$, it follows from a result of
Tumanov [Tu] that for every compact set $K_l \subset U_l$
both
$F$ and
$G$ extend holomorphically into an open connected wedge
$\Omega_l$ in $\bC^K$ with edge on $K_l$. Also,  since
$M^\prime$ is CR diffeomorphic with $M$,  it follows that
the intrinsic complexification 
$\scrV^\prime$ of
$M^\prime$ has complex dimension $K$ as well. Since $M$
is 
$l(M)$-nondegenerate at $p_l$ and hence also essentially
finite at $p_l$ (see [BER1]), it follows that both $F$ and $G$
extend as biholomorphisms of some neighborhood of 
$p_l$ in $\bC^K$ onto a neighborhood of $F(p_l)\in
M^\prime$ in $\scrV^\prime$ (see [BJT]). It follows from
\thetag{2.1} and  Proposition 2.3 that the holomorphic
extensions of $F$ and $G$ into $\Omega_l$ are identical.
Consequently, $F\equiv G$ in $U_l$ and, since this is
true for any 
$l=1,...,k$, the theorem follows by continuity of the
mappings. \qed
\enddemo

\heading 3. The infinitesimal CR automorphisms\endheading

\subhead 3.1. The minimal case\endsubhead  We shall prove
Theorem 2 in this section. 
\comment
\noindent{\bf Remarks:} \roster
\item"(i)" If we, as in the remark (ii) following Theorem
2.1.1,  denote by $k$ the number of components of
$M\setminus V$, where $V$ is the proper real analytic
subset of $M$ consisting of those points at which $M$ is
not minimal, that have $p_0$ in their closure then it
follows from the proof that the dimension of 
$\aut(M,p_0)$ is bounded by $k\cdot\dim M$ times the
number of  monomials in $\dim_\bR M$ variables of degree
$\leq Nl(M)$.
\item"(ii)" If $M$ is not minimal anywhere,
$p_0\in M$ is a point at which the local CR orbit has
maximal dimension, and
$\aut(M,p_0)\neq\{0\}$ then \thetag{3.1.2} also holds as
will be apparent from the proof. \endroster

\medskip

\demo{Proof} Let us begin by proving that $\hol(M,p_0)$
is infinite dimensional whenever $M$ is holomorphically
degenerate, or everywhere non-minimal and  homogeneous.
The proof of this in the case where $M$ is
holomorphically  degenerate is exactly the same as in the
hypersurface case (see [St2]).  If $M$ is homogeneous and
everywhere non-minimal then there is a  holomorphic
polynomial $h(z,w)$ which is real valued and non-constant
on $M$  (see [BER1]). Also, there is a non-trivial vector
field $X\in\hol(M,p_0)$ (e.g. the infinitesimal dilation;
see [St1]). The infinite dimensionality of $\hol(M,p_0)$
follows by noting  that $h^kX\in\hol(M,p_0)$, for
$k=0,1,...$, are all linearly independent over $\bR$. Let
us now prove \thetag{3.1.1} under the assumptions in the
first part of the theorem.
\endcomment
Let $X^1,...,X^m$ $\in\aut(M,p_0)$ be linearly independent
over $\bR$. Let
$x=(x_1,...,x_r)$ be a local coordinate system on $M$
vanishing at $p_0$. In this coordinate system, we may
write
$$ X^j=\sum_{l=1}^r\tilde
X^j_l(x)\frac{\partial}{\partial x_l}=
\tilde X^j(x)\cdot\frac {\partial}{\partial x}.\tag3.1.1
$$ For
$y=(y_1,...,y_m)\in\bR^m$, we denote by $\Phi(t,x,y)$ the
flow of the vector field $y_1X_1+...+y_mX_m$, i.e. the
solution of
$$
\left\{\aligned&\frac{\partial\Phi}{\partial t}(t,x,y)=
\sum_{i=1}^my_i\tilde
X^i(\Phi(t,x,y))\\&\Phi(0,x,y)=x.\endaligned\right.
\tag3.1.2
$$ By choosing $\delta>0$ sufficiently small, there is
$c>0$ such that the flows $\Phi(t,x,y)$ are smooth
($C^\infty$) in 
$\{(t,x,y)\in\bR^{1+r+m}\:|t|\leq 2,\,|x|\leq
c,\,|y|\leq\delta\}$. This follows from the identity
$$
\Phi(st,x,y)=\Phi(t,x,sy),\quad s\in\bR,\tag3.1.3
$$ which, in turn, follows from the fact that the
solution of \thetag{3.1.2} is unique (the reader can
verify that the left side of \thetag{3.1.3} solves the
initial value problem \thetag{3.1.2} that defines the
right side of \thetag{3.1.3}). Denote by $F(x,y)$ the
corresponding time-one maps, i.e.
$$ F(x,y)=\Phi(1,x,y).\tag3.1.4
$$
\proclaim{Assertion 3.1.1} There is a $\delta^\prime$, 
$0<\delta^\prime<\delta$, such that for any fixed
$y_1,y_2$ with $|y^1|,|y^2|\leq \delta^\prime$, if
$F(x,y^1)\equiv F(x,y^2)$ for
$|x|\leq c$ then
necessarily
$y^1=y^2$.\endproclaim
\demo{Proof of Assertion 3.1.1} Note from \thetag{3.1.3}
that, with the notation $e_i$ for $i$th unit vector in
$\bR^m$, 
$$
\frac{\partial F}{\partial y_i}(x,0)  = \frac{d}{d
s}[\Phi(1,x,se_i)]_{s=0}
=\frac{\partial\Phi}{\partial t}(0,x,e_i)=\tilde X^i(x).
\tag 3.1.5
$$ 
Thus, denoting by $\tilde X(x)$ the $r\times m$-matrix
with column vectors
$\tilde X^i(x)$, we have
$$
\frac{\partial F}{\partial y}(x,0)=\tilde X(x).\tag3.1.6
$$  If $J(x)=(J_1(x),...,J_r(x))$ is a continuous mapping
into $\bR^r$ for $|x|\leq c$ then we put
$$ ||J(x)||=\sup_{|x|\leq
c}\left(\sum_{l=1}^rJ_l(x)^2\right)^{1/2}.\tag3.1.7
$$ By Taylor expansion  we obtain
$$ ||F(x,y^2)-F(x,y^1)||\geq\left|\left|\frac{\partial
F}{\partial y}(x,y^1)
\cdot (y^2-y^1)\right|\right|-C|y^2-y^1|^2,\tag3.1.8
$$ where $C>0$ is some uniform constant for
$|y^1|,|y^2|\leq\delta$.  Now, by assumption, the vector
fields
$X^1,...,X^m$ are linearly independent over $\bR$. This
means precisely that there is a constant $C^\prime$ such
that
$$ ||\tilde X(x)\cdot y||\geq C^\prime|y|.\tag3.1.9
$$ Using 
\thetag{3.1.6}, \thetag{3.1.9}, the smoothness of $F$,  
and a standard compactness argument, we deduce from
\thetag{3.1.8} that the  conclusion of Assertion 3.1.1
holds.\qed

Now, we proceed with the proof of Theorem 2. Denote by
$U$ the open neighborhood of $p$ on $M$ given by
$|x|<c$.  We make use of Theorem 2.1.1 with $M$ replaced
by
$U$.  Let
$p_1,...,p_k$ be the points in
$U$ given by the theorem. By choosing the number
$\delta^\prime>0$ in Assertion 3.1.1 even smaller if
necessary, we may assume that the maps
$x\mapsto F(x,y)$, for $|y|<\delta^\prime$, are CR
diffeomorphisms of $U$ into $M$. Consider the smooth
mapping from $|y|<\delta^\prime$ into $\bR^\mu$ defined by
$$ y\mapsto\left(\frac{\partial^{|\alpha|}
F(p_l,y)}{\partial x^\alpha}
\right)\in \bR^\mu,\tag3.1.12
$$ where $l=1,...,k$, and $|\alpha|\leq (d+1)l(M)$ (thus,
the dimension $\mu$ equals $k\cdot r$ times the
number of monomials in $r$ variables of degree 
$\leq (d+1)l(M)$). This mapping is injective for
$|y|<\delta^\prime$ in view of Theorem 2.1.1 and
Assertion 3.1.1. Consequently, we have a smooth injective
mapping from a neighborhood of the origin in $\bR^m$ into
$\bR^\mu$. This implies that $m\leq \mu$ and hence the
desired finite dimensionality of the conclusion of
Theorem 2.
\qed\enddemo

\subhead 3.2. The non-minimal case\endsubhead  We shall
prove a  generalization of Theorem 2 for the case where
$M$ is not minimal.  Before stating this result, we need
some notation. Let
$M\subset\bC^N$ be a real analytic CR submanifold, and
let $p_0\in M$. Denote by $\Cal F_M(p_0)$ the ring of
germs of $C^\infty$ real valued functions on $M$ at $p_0$
which are also CR,  and by $\Cal G_M(p_0)$ the subring
consisting of those that are real analytic, i.e. those
that are restrictions of the elements in $\scrO_N(p_0)$
which  are real valued on $M$. If $f$ is a representative
of an element in 
$\Cal F_M(p_0)$ then the restriction of $f$ to each CR
orbit is constant. (Conversely, if $f$ is a real-valued
$C^\infty$ function near $p_0\in M$ which  is constant on
the CR orbits then $f\in\Cal F_M(p_0)$.) Hence, if
$M$ is connected and minimal somewhere then $\Cal
G_M(p_0)=
\Cal F_M(p_0)=\bR$. On the other hand, if $M$ is
non-minimal everywhere then it follows from the Frobenius
theorem that $\Cal G_M(p)$ contains a non-trivial element
if the CR orbit of $p$ has maximal dimension. 

It is also easy to verify that $\aut(M,p)$ is a $\Cal
F_M(p)$-module, and $\hol(M,p)$ is a $\Cal
G_M(p)$-module, for every $p\in M$. Indeed, it is obvious
that
$\hol(M,p)$ is a $\Cal G_M(p)$-module, and the fact that
$\aut(M,p)$ is a $\Cal F_M(p)$-module is also immediate
from the following characterization  of $\aut(M,p)$ (see
e.g. [BR2]): {\it a smooth real vector field $X$ on $M$ near
$p\in M$ is in $\aut(M,p)$ if and only if}
$$ [\bL,X]\in\bL,\tag3.2.1
$$ {\it where $\bL$ denotes the space of smooth CR vector
fields on $M$ near $p$.}

The following is a generalization of Theorem 2.

\proclaim{Theorem 3.2.1} Let $M\subset\bC^N$ be a real
analytic,  holomorphically nondegenerate CR submanifold,
and assume that $M$ is  everywhere non-minimal. Then
there is a dense open subset $\Omega_a\subset M$ such
that $\aut(M,p)$ is a finitely generated free
$\Cal F_M(p)$-module for $p\in \Omega_a$, and a dense
open subset 
$\Omega_h\subset M$ such that $\hol(M,p)$ is a finitely
generated free
$\Cal G_M(p)$-module for $p\in \Omega_h$.\endproclaim We
begin with a local description of generic manifolds which
are everywhere non-minimal.
\proclaim{Proposition 3.2.2} Let $M\subset\bC^N$ be a
generic, real analytic submanifold of codimension $d$ which
is everywhere non-minimal, and let $p_0\in M$ whose local CR
orbit is of maximal dimension. Then there are coordinates
$(z,w^\prime,
w^{\prime\prime})\in\bC^n\times\bC^{d-q}\times\bC^q=\bC^N$,
where $q$ denotes the codimension of the local CR orbit of
$p_0$ in $M$, vanishing at $p_0$  such that $M$ is
defined by the equations
$$
\left\{\aligned&\im w^\prime=\phi(z,\bar z,\re
w^\prime,\re w^{\prime\prime})
\\&\im w^{\prime\prime}=0;\endaligned\right.\tag3.2.2
$$ here, $\phi$ is a real valued analytic function with
$\phi(z,0,s^\prime, s^{\prime\prime})\equiv0$. Moreover,
the local CR orbit of the point
$(z,w^\prime,w^{\prime\prime})=(0,0,s'')$, for
$s''\in\bR^q$, is given by 
$$
\left\{\aligned&\im w^\prime=\phi(z,\bar z,\re
w^\prime,s'')
\\&w^{\prime\prime}=s'',\endaligned\right.\tag3.2.3
$$
\endproclaim
\demo{Proof} Since the local CR orbit of $p_0\in M$ has
maximal dimension, it follows from the Frobenius theorem
that there are $h_1,...,h_q\in
\scrO_N(p_0)$ such that 
$$
\partial h_1\wedge...\wedge \partial h_q\neq 0\tag3.2.4
$$ near $p_0$, the restriction of each $h_j$ to $M$ is
real valued, and such that the CR orbit of any point
$p_1$ near $p_0$ is given by $\{p\in
M\:h_1(p)=h_1(p_1),...,h_q(p)=h_q(p_1)\}$. Thus, there
are coordinates
$(Z^\prime,w^{\prime\prime})\in\bC^{n+d-q}\times
\bC^q$, vanishing at $p_0$, such $M$ is contained in 
$\im w^{\prime\prime}=0$ and such that the CR orbits on
$M$ are given as the intersections between 
$M$ and $w^{\prime\prime}=s''$, for $s''$ near
$0\in\bR^q$. Since $M$ is generic and contained in the
flat surface
$\im w^{\prime\prime}=0$,  we can make a change of
coordinates $\tilde Z^\prime=A(w^{\prime\prime})
Z^\prime$, where $A(w^{\prime\prime})$ is an invertible 
$(n+d-q)\times(n+d-q)$-matrix with holomorphic matrix
elements in 
$w^{\prime\prime}$, and write $\tilde
Z^\prime=(u,v)\in\bC^n\times
\bC^{d-q}$ such that $M$ is given by
$$
\left\{\aligned&\im v=\psi(u,\bar u,\re v,\re
w^{\prime\prime})
\\&\im w^{\prime\prime}=0,\endaligned\right.\tag3.2.5
$$ where $\psi(u,\bar u,s^\prime,s^{\prime\prime})$ is
real valued, real analytic, and satisfies 
$$
\psi(0,0,0,s^{\prime\prime})\equiv\frac{\partial\psi}{\partial
u} (0,0,0,s^{\prime\prime})\equiv
\frac{\partial\psi}{\partial s^\prime}
(0,0,0,s^{\prime\prime}) \equiv 0\tag3.2.6
$$ for all $s^{\prime\prime}\in\bR^q$ near $0$. Now, we claim
that we can actually find  holomorphic coordinates
$(z,w^\prime)\in\bC^n\times\bC^{d-q}$ of the form
$$
\left\{\aligned&z=u\\&w^\prime=f(u,v,w^{\prime\prime})
\endaligned\right.\tag3.2.7
$$ such that $M$ is defined by
$$
\left\{\aligned&\im w^\prime=\phi(z,\bar z,\re
w^\prime,\re w^{\prime\prime})
\\&\im w^{\prime\prime}=0,\endaligned\right.\tag3.2.8
$$ where $\phi(z,\bar z,s^\prime,s^{\prime\prime})$ is
real valued, real analytic, and satisfies 
$$
\phi(z,0,s^\prime,s^{\prime\prime})\equiv0.\tag3.2.9
$$ This follows from [BJT, Lemma 1.1] in the following
way. Consider only the equation for $\im v$ in
\thetag{3.2.5} for a fixed $s^{\prime\prime}=
\re w^{\prime\prime}$.  Lemma 1.1 in [BJT] asserts that
there is a change of coordinates
$$
\left\{\aligned&z=u\\&w^\prime=f(u,v;s^{\prime\prime}),
\endaligned\right.\tag3.2.10
$$ holomorphic in $(u,v)$, such that the equation for
$\im v$ in \thetag{3.2.5} becomes the equation for $\im
w^{\prime}$ in \thetag{3.2.8} and $\phi$ is as
in \thetag{3.2.9}. Moreover, the change of
coordinates is  obtained by the implicit function
theorem, so if we now think of $s^ {\prime\prime}$ as a
real analytic parameter in \thetag{3.2.5} we find that
the change of coordinates \thetag{3.2.10} is real
analytic in 
$s^{\prime\prime}$.  Hence, it can be extended as a
holomorphic change of coordinates of the form
\thetag{3.2.7}. This change of coordinates coincides with 
\thetag{3.2.10} on $\im w^{\prime\prime}=0$, and since
$M$ is contained in $\im w^{\prime\prime}=0$ the claim is
proved. This completes the proof of the
proposition.\qed\enddemo

The following is an easy consequence of Proposition 3.2.2
and the definition of $k$-nondegeneracy.

\proclaim{Corollary 3.2.3} Let $M\subset\bC^N$ be a real
analytic, holomorphically nondegenerate CR submanifold
and $p \in M$.  If the local CR orbit
$W_p$ of $p$ is of maximal dimension then, for any
integer $k$,
$M$ is
$k$-nondegenerate at $p$ if and only if $W_p$ is
$k$-nondegenerate at $p$. In particular, $l(M)=l(W_p)$. 
\endproclaim

\demo{Proof of Theorem 3.2.1} We start by showing that for
$p$ in an open dense subset $\Omega_a \subset M$,
$\aut(M,p)$ is  a finitely generated free $\Cal
F_M(p)$-module. Denote by
$W_p$ the local CR orbit of any point $p\in M$. Let
$\Omega^1_a$ be the dense open subset of points $p\in M$ such
that the dimension of $W_p$ is maximal and such that $M$
is $l(M)$-nondegenerate at $p$. By Corollary  3.2.3, the
CR orbit $W_p$, for $p\in\Omega_a^1$, is also
$l(M)$-nondegenerate at $p$ and $l(W_p) = l(M)$. Next,
define the following function on $\Omega^1_a$
$$
\alpha(p) = \dim_\bR
\left(\aut(M,p)|_p/\tilde T_p(W_p)\right),\tag3.2.11
$$ where $\aut(M,p)|_p$ denotes the subspace of $T_p(M)$
obtained as values of the vector fields in 
$\aut(M,p)$ and $\tilde T_p(W_p) = \aut(M,p)|_p\cap
T_p(W_p)$.  Let $\Omega^2_a$ be the  set of points
$p\in
\Omega^1_a$ such that both $\dim_\bR \aut(M,p)|_p$ and 
$\alpha(p)$ are maximal in a neighborhood of
$p$.  Using the fact that the the CR orbits $W_q$ for $q$
in a neighborhood of $p\in \Omega^1_a$ form a real
analytic foliation  of $M$, one can check by elementary
linear algebra that $\Omega^2_a$ is open and dense in
$\Omega^1_a$ and that $\alpha(q)$ is constant for $q$ in a
neighborhood of $p$. 

Let us denote
by $\aut(M,p)^\prime$ the subspace of
$\aut(M,p)$ consisting of those vector fields that are
also tangent to $W_q$ for all $q$ in a neighborhood of
$p$. It follows immediately from \thetag{3.2.1} that the
restriction of any $X\in \aut(M,p)^\prime$ to $W_q$, for
$q$ in a  neighborhood of $p$, is in $\aut(W_q,q)$. We 
define the following function in $\Omega^1_a$
$$
\beta(p)=\dim_\bR \aut(M,p)^\prime|_{W_p},\tag3.2.12
$$ where $\aut(M,p)^\prime|_{W_p}$ denotes the subspace
of $\aut(W_p,p)$ obtained by taking restrictions to $W_p$
of the vector fields in 
$\aut(M,p)^\prime$. Let $d$ denote the codimension of the
maximal CR orbits in $\bC^N$. It follows from Theorem 1 and
Corollary 3.2.3 that
$0\leq
\beta(p)<\mu$, where $\mu$ denotes the dim$_\bR M$ times the
number of monomials in
$\dim_\bR M$ variables of degree $\leq (d+1)l(M)$, for all 
$p\in \Omega^1_a$. Clearly, if we have vector fields
$Y_1,...,Y_k$ on 
$M$ that are tangent to each $W_q$, for $q$ near
$p\in\Omega^1_a$, and linearly independent over 
$\bR$ as  vector fields  on $W_p$ then they are also
linearly independent over $\bR$ on $W_q$, for all $q$
in a neighborhood of $p$.  We  let
$\Omega^3_a$ be the set of points $p\in\Omega^1_a$ for
which $\beta(p)$ is maximal in a neighborhood of
$p$. It follows easily that $\Omega^3_a$ is a dense open
subset of
$\Omega^1_a$ and that $\beta(q)$ is constant for $q$ in
a neighborhood of $p$.

Now, let $\Omega_a$ be the intersection of $\Omega^2_a$
and $\Omega^3_a$ in
$\Omega^1_a$. This is a dense open subset of $M$, and we
claim that 
$\aut(M,p_0)$ is a finitely generated free $\Cal
F_M(p_0)$-module for every
$p_0\in\Omega_a$. Let
$X^\prime_1,...,X^\prime_{\beta(p_0)},T_1,...,T_
{\alpha(p_0)}
\in\aut(M,p_0)$ be such that the images of
$X^\prime_1,...,X^\prime_ {\beta(p_0)}$ in
$\aut(M,p_0)^\prime|_{W_{p_0}}$ form a basis for that
space, and the images of
$T_1,...,T_{\alpha(p_0)}$ in
$\aut(M,p_0)|_{p_0}/\tilde T_{p_0}(W_{p_0})$ form a  basis
for the latter. As we noted above, the images of these
vector fields  are linearly independent in the
corresponding vector spaces at every $p$ in a
neighborhood of $p_0$. This  implies that these vector
fields are linearly independent over the ring
$\Cal F_M(p_0)$. To see this, assume that there are
$c_1,...,c_{\alpha(p_0)},d_1,...,d_{\beta(p_0)}\in\Cal
F_M(p_0)$ such that
$$
\sum_{i=1}^{\beta(p_0)}d_iX^\prime_i+\sum_{j=1}^{\alpha(p_0)}c_jT_j\equiv0.
\tag3.2.13
$$ Taking the image of \thetag{3.2.13} in
$\aut(M,p)|_p/\tilde T_p(W_p)$, for $p$ near
$p_0$, and using the fact that the $X^\prime_i$ are
tangent to $W_p$, we  deduce that $c_1\equiv...\equiv
c_{\alpha(p_0)}\equiv0$. Then taking the image in
$\aut(M,p)|_{W_p}$ and using the fact that the $d_i$ are
constant on $W_p$, we deduce that $d_1\equiv...\equiv
d_{\beta(p_0)}\equiv0$. Hence, the vector fields are
linearly independent over $\Cal F_M(p_0)$.

It remains to prove that these vector fields generate
$\aut(M,p_0)$ as a $\Cal F_M(p_0)$-module. Now, since
both $\alpha(p)$ and
$\beta(p)$ are constant in a neighborhood of $p_0$, it
follows that the images of
$X^\prime_1,...,X^\prime_{\beta(p_0)}$ in 
$\aut(M,p)^\prime|_{W_p}$ form a basis for this vector
space for all $p$ in a neighborhood of $p_0$, and
similarly for the images of $T_1,\dots,T_ {\alpha(p_0)}$.
Hence, the vector spaces $\aut(M,p)^\prime|_{W_p}$ and
$\aut(M,p)|_p/\tilde T_p(W_p)$ form $C^\infty$ real vector
bundles (not real analytic since the coefficients of
these vector fields are merely
$C^\infty$ functions) over $M$ in a neighborhood of
$p_0$. Thus, if we take any $X\in\aut(M,p_0)$ and
consider its image in the latter vector bundle we obtain
smooth real valued functions $c_1,...,c_{\alpha(p_0)}$ on
$M$ such that
$$ X^\prime=X-\sum_{j=1}^{\alpha(p_0)}c_j T_j\tag3.2.14
$$ is tangent to $W_p$ for all $p$ in a neighborhood of
$p_0$. We claim that 
$$
\align& X^\prime\in\aut(M,p_0)^\prime\tag3.2.15
\\&c_1,...,c_{\alpha(p_0)}\in\Cal
F_M(p_0).\tag3.2.16\endalign
$$ To verify this, we compute the bracket of $X$ with a
CR vector field $L\in\bL$
$$
\aligned [L,X]
&=[L,X^\prime]+\sum_{j=1}^{\alpha(p_0)}[L,c_jT_j]\\
&=[L,X^\prime]+\sum_{j=1}^{\alpha(p_0)}c_j[L,T_j]+(Lc_j)T_j.
\endaligned\tag3.2.17
$$ Since $X,T_1,...,T_{\alpha(p_0)}\in\aut(M,p)$, we know
that $[L,X],[L,T_1],..., [L,T_{\alpha(p_0)}]\in\bL$.
Also, since both $X^\prime$ and $L$ are tangent to $W_p$
for all $p$ near $p_0$, \thetag{3.2.17} implies that
$$
\sum_{j=1}^{\alpha(p_0)}(Lc_j)T_j\tag3.2.18
$$ is tangent to $W_p$ for all $p$ near $p_0$. Since the
images of $T_1,..., T_{\alpha(p_0)}$ in the vector space
$\aut(M,p)|_p/\tilde T_p(W_p)$ are linearly  independent
for
$p$ near $p_0$, we deduce  that
$Lc_1(p),...,Lc_{\alpha(p_0)}(p)=0$ for all $p$ in a
neighborhood of $p_0$. Hence, \thetag{3.2.16} is proved,
because
$L\in \bL$ was arbitrary. Moreover, \thetag{3.2.15} also
follows because now \thetag{3.2.1}  implies that
$[L,X^\prime]\in\bL$. The claim is proved. 

To finish the proof of the first part of the theorem, we
have to prove that $X^\prime\in
\aut(M,p_0)^\prime$ can be written
$$ X^\prime=\sum_{i=1}^{\beta(p_0)}d_j X_j,\tag3.2.19$$
where $d_1,...,d_{\beta(p_0)}\in\Cal F_M(p_0)$. By taking
the image in the real $C^\infty$ vector bundle
$\aut(M,p)^\prime|_{W_p}$, we obtain smooth  real valued
functions $d_1,...,d_{\beta(p_0)}$ such that
\thetag{3.2.19} holds. Since the values
$d_1,...,d_{\beta(p_0)}$ are unique as real numbers on
each CR orbit, it follows that each function $d_j(p)$ is
constant on the CR orbits near $p_0$ and, hence,
$d_j\in\Cal F_M(p_0)$. This completes the proof of the
statement that
$\aut(M,p_0)$ is a finitely generated free $\Cal
F_M(p_0)$-module at every $p_0\in\Omega_a$. 

To prove the corresponding statement for $\hol (M,p)$ we
define
$\Omega_h$ in complete analogy with $\Omega_a$, replacing
$\aut(M,p)$ by $\hol(M,p)$. The same proof as above, {\it
mutatis mutandi},  completes the proof of the statement
that $\hol(M,p_0)$ is a finitely generated free $\Cal
G_M(p_0)$-module, since the vector bundles corresponding
to those with fibers $\aut(M,p)|_p/T_p(W_p)$ and
$\aut(M,p)^\prime|_{W_p}$, replacing $\aut(M,p)$ by
$\hol(M,p)$, are real analytic. We leave the details of
this to the reader.\qed\enddemo

\head 4 Proof of Theorem 3 and examples \endhead

\subhead 4.1 Proof of Theorem 3 \endsubhead
 To prove Theorem 3, suppose first that
$M$ is holomorphically degenerate. Then for any
$p\in M$ there is a  germ of a vector field
$X=\sum_{j=1}^N c_j(Z) {\partial \over
\partial Z_j}$, with $c_j(Z)$ holomorphic, tangent to
$M$ near $p$ with nonzero restriction to $M$. Let $h(Z)$ be
any holomorphic function nontrivial on $M$ with $h(p) =0$ and
$K $ a positive integer as in the statement of the theorem. 
Then the vector field $Y = \re (h(Z)^K X)$ is a nontrivial
element of
$\hol_\bR(M,p_0)$ and vanishes to order at least $K$ at
$p$.  For
$t>0$ small the CR mapping $F=\exp tY$ extends to a local
biholomorphism mapping $M$ into itself with $F(p)=p$. It is
not hard to check that all derivatives of $F$ up
to order
$K$ agree with those of the identity mapping, $I$, yet
$F\not\equiv I$ on $M$.  This proves (0.3) with $G = I$.

Next, we show that $\hol_\bR(M,p)$ is infinite
dimensional whenever
$M$ is holomorphically degenerate, or everywhere non-minimal
and  homogeneous.  The proof of this in the
case where
$M$ is holomorphically  degenerate is exactly the same as in
the hypersurface case (see [St2]).  If $M$ is homogeneous and
everywhere non-minimal then there is a  holomorphic
polynomial $h(z,w)$ which is real valued and non-constant on
$M$  (see [BER1]). Also, there is a non-trivial vector field
$X\in\hol(M,p)$ (e.g. the infinitesimal dilation; see
[St1]). The infinite dimensionality of $\hol(M,p)$ follows
by noting  that $h^kX\in\hol(M,p_0)$, for $k=0,1,...$, are
all linearly independent over $\bR$.  We note also that by
considering the mappings $\exp( t\ h(z,w)^K X )$, we can
construct nontrivial local biholomorphisms which agree
with the identity up to any preassigned order. 

Finally, we shall show that if $M$ is not minimal almost
everywhere then for $p\in M$ outside an open dense subset
 if {\rm dim}$_\bR \hol(M,p)
\not= 0$, then {\rm dim}$_\bR\hol(M,p)
= \infty$. Indeed, suppose the local CR orbit of $p$ is of
maximal dimension. Then it follows from Proposition
3.2.2 that there is a  holomorphic function $h(Z)$
whose restriction to $M$ is real and nonconstant. If $X$
is a nontrivial vector field in $\hol(M,p)$, then so is
$X_k = h(Z)^k X $ for any positive integer $k$. 
Since the $X_k$ are all linearly independent as vector
fields over $\bR$, it follows that {\rm dim}$_\bR\hol(M,p)
= \infty$.  This completes the proof of Theorem 3.\qed

\subhead 4.2. A holomorphically nondegenerate, nowhere
minimal CR submanifold with $\hol(M,0)=\{0\}$\endsubhead
 In Theorem 3, it is shown that for most points $p \in M$,
if $\hol (M,p)$ contains at least one
non-trivial element and if $M$ is nowhere minimal, then
{\rm dim}$_\bR\hol(M,p) = \infty$. In this section we 
construct an example to show 
it may  happen that
$\hol(M,p)=\{0\}$ for a holomorphically nondegenerate,
nowhere minimal CR submanifold $M$.

N. Stanton [St3] has given examples of real
hypersurfaces with no (non-trivial) infinitesimal CR
automorphisms. We leave it to the reader  to verify that
a slight modification of the argument of [St3]
proves the following.
\proclaim{Proposition 4.2.1 ([St3])} Let
$M^0\subset\bC^2$ be the  hypersurface defined by
$$
\im w=z^4\bar z^{10}+z^{10}\bar z^4+(\re w)|z|^8.\tag4.2.1
$$ Then $\aut(M^0,0)=\{0\}$ (and hence also
$\hol(M,0)=\{0\}$).\endproclaim

We will use this to prove:
\proclaim{Proposition 4.2.2} Let $M\subset\bC^3$ be
defined by
$$
\left\{\aligned&\im w_1=z^4\bar z^{10}+z^{10}\bar
z^4+(\re w_1)|z|^8+ (\re w_2)|z|^4\\&\im
w_2=0.\endaligned\right.\tag4.2.2
$$ Then $\hol (M,0)=\{0\}$.\endproclaim
\noindent{\bf Remark:} Note that $M$ is holomorphically
nondegenerate and nowhere minimal.\medskip

\demo{Proof} Assume that
$$ X=a(z,w)\frac{\partial}{\partial z}+\bar a(\bar z,\bar
w)\frac{\partial} {\partial \bar
z}+b(z,w)\frac{\partial}{\partial w_1}+
\bar b(\bar z,\bar w)\frac{\partial}{\partial \bar w_1}+
c(z,w)\frac{\partial}{\partial w_2}+
\bar c(\bar z,\bar w)\frac{\partial}{\partial \bar
w_2}\tag4.2.3
$$ is in $\hol (M,0)$. This is the same as saying that
the holomorphic vector field
$$ Y=a(z,w)\frac{\partial}{\partial z}+\bar
a(\chi,\tau)\frac{\partial} {\partial
\chi}+b(z,w)\frac{\partial}{\partial w_1}+
\bar b(\chi,\tau)\frac{\partial}{\partial \tau_1}+
c(z,w)\frac{\partial}{\partial w_2}+
\bar c(\chi,\tau)\frac{\partial}{\partial
\tau_2}\tag4.2.4
$$ in $\bC^6$ is tangent to the complexification $\scrM$
of $M$ in $\bC^6$
$$
\left\{\aligned&w_1-\tau_1-2i(z^4\chi^{10}+z^{10}\chi^4)-i(w_1+\tau_1)z^4
\chi^4-i(w_2+\tau_2)z^2\chi^2=0\\&
w_2-\tau_2=0.\endaligned\right.\tag4.2.5
$$ Let us denote the first equation above by
$\rho(z,w,\chi,\tau)=0$. It is an easy exercise to verify
that if $Y$ is tangent to $\scrM$ then $b(z,w)$ and
$c(z,w)$ are independent of $z$ and real, i.e.
$b=\bar b$ and $c=\bar c$. We will use the notation
$b=b(w)$ and 
$c=c(w)$. We let
$$ Y^\prime=a(z,w)\frac{\partial}{\partial z}+\bar
a(\chi,\tau)\frac{\partial} {\partial
\chi}+b(w)\frac{\partial}{\partial w_1}+
\bar b(\tau)\frac{\partial}{\partial \tau_1}\tag4.2.6
$$ and 
$$
\rho_0(z,w_1,\chi,\tau_1)=w_1-\tau_1-2i(z^4\chi^{10}+z^{10}\chi^4)-i(w_1+
\tau_1)z^4\chi^4.\tag4.2.7
$$ Note that $\rho_0=0$ is the defining equation of the
complexification
$\scrM^0$ of $M^0$ in $\bC^4$. Applying $Y$ to $\rho$, we
obtain
$$
\aligned(Y\rho)(z,w,\chi,\tau)=&
(Y^\prime\rho_0)(z,w,\chi,\tau)- ia(z,w)(w_2+\tau_2)
z\chi^2\\&-i\bar
a(\chi,\tau)(w_2+\tau_2)z^2\chi-i(c(w)+c(\tau))z^2\chi^2.
\endaligned
\tag4.2.8
$$ Since $Y$ is assumed tangent to $\scrM$, this
expression is 0 on $\scrM$. We can solve for $\tau$ in
the defining equations of $\scrM$ and obtain
$$
\tau_1=\bar
Q_0(\chi,z,w_1)+O(w_2)\quad,\quad\tau_2=w_2,\tag4.2.9
$$ where $\tau_1=\bar Q_0(\chi,z,w_1)$ is the defining
equation of $\scrM^0$ in $\bC^4$ and $O(w_2)$ denotes, as
usual, terms that contain the factor
$w_2$. We have also substituted $\tau_2=w_2$ in the first
equation. Let us write \thetag{4.2.9} as
$\tau=\bar Q(\chi,z,w)$ for short. Note that $\bar
Q(0,z,w)\equiv\bar Q(\chi,0,w)\equiv w$ and $\bar
Q_0(0,z,w_1)\equiv\bar Q_0(\chi,0,w_1)\equiv w_1$.
Substituting in \thetag{4.2.8} we obtain
$$
\aligned&(Y^\prime\rho_0)(z,w,\chi,\bar
Q(\chi,z,w))\equiv\\& 2ia(z,w)w_2z\chi^2+ 2i\bar
a(\chi,\bar Q(\chi,z,w))w_2z^2\chi+i(c(w)+c(\bar
Q(\chi,z,w)))z^2
\chi^2.\endaligned\tag4.2.10
$$ Let us expand the holomorphic vector field $Y^\prime
(z,w_1,w_2,\chi,\tau_1,w_2)$ in $w_2$. We obtain
$$ Y^\prime(z,w_1,w_2,\chi,\tau_1,w_2)=\sum_{k=k_0}^\infty
Y^\prime_k(z,w_1,\chi,\tau_1)w_2^k,\tag4.2.11
$$ where each $Y^\prime_k$ is a holomorphic vector field
in the variables $(z,w_1,\chi,\tau_1)$ and where in
particular the vector field $Y^\prime_{k_0}$ is not
identically 0; we assume here, in  order to obtain a
contradiction, that $Y^\prime$ is not identically 0. Note
that, since $k_0$ is assumed to be the lowest order in
the expansion of
$Y^\prime(z,w_1,w_2,\chi,\tau_1,w_2)$, the coefficients
$a(z,w_1,w_2)$ and $b(w_1,w_2)$ have to be divisible by
$w_2^{k_0}$ ($k_0$ could, of course, be 0). We expand
$a,b,c$ in $w_2$ as follows
$$ a(z,w)=\sum_{k=k_0}^\infty
a_k(z,w_1)w_2^k\,,\,b(w)=\sum_{k=k_0}^\infty
b_k(w_1)w_2^k\,,\,c(w)=\sum_{k=0}^\infty
c_k(w_1)w_2^k.\tag4.2.12
$$ Identifying the coefficients of the lowest order term
in $w_2$ (i.e. of $w_2^{k_0}$) in \thetag{4.2.10}, using
the fact that
$$
\bar Q(\chi,z,w)=(\bar Q_0(\chi,z,w_1)+O(w_2),w_2),
$$  we find 
$$ (Y^\prime_{k_0}\rho_0)(z,w_1,\chi,\bar
Q_0(\chi,z,w_1))=i(c_{k_0}(w_1)+
c_{k_0}(Q_0(\chi,z,w_1)))z^2\chi^2.\tag4.2.13
$$ Now, it is easy to see that the expansion of the left
hand side in terms of $z$ and $\chi$ does not contain a
term with $z^2\chi^2$. The expansion of the right hand
side contains the term
$$ 2ic_{k_0}(w_1)z^2\chi^2.\tag4.2.14
$$ Thus, we must have $c_{k_0}(w_1)\equiv0$ and
$$ (Y^\prime_{k_0}\rho_0)(z,w_1,\chi,\bar
Q_0(\chi,z,w_1))\equiv0.\tag4.2.15
$$ The latter implies that the vector field
$Y^\prime_{k_0}$ in $\bC^4$ is tangent to $\scrM^0$ or,
equivalently, that the vector field
$X^\prime_{k_0}$ in $\bC^2$, obtained by formally
replacing $\chi$ by $\bar z$ and $\tau_1$ by $\bar w_1$
in $Y^\prime_{k_0}$, is tangent to $M^0$. Now, 
$X^\prime_{k_0}$ is the real part of a holomorphic vector
field so, since $X^\prime_{k_0}$ is tangent to $M^0$,
$X^\prime_{k_0}\in\hol (M^0,0)$ and hence
$X^\prime_{k_0}\equiv0$ by Proposition 4.2.1. This
contradicts the fact that $Y^\prime_{k_0}$ was assumed
$\not\equiv0$. Consequently,
$Y^\prime$ is identically 0. That means $Y$ has to be of
the form
$$ c(w)\frac{\partial}{\partial
w_2}+c(\tau)\frac{\partial}{\partial \tau_2}.
\tag4.2.16
$$ It is easy to check that this implies $c(w)\equiv0$ as
well. This completes the proof of Proposition .\qed\enddemo

\head \S 5 Remarks\endhead

We shall restrict our remarks to the case where $M$ is
a generic manifold. For any $p \in M$ we let $G_p$ denote
the set of germs $H$ of  biholomorphisms near $p$, with
$H(M) \subset M$ and $H(p) = p$.  It is easy to see that
the set $G_p$ forms a group under composition of mappings. 
We have the following.  

\proclaim {Theorem 5.1}
Let $M$ be a real analytic, holomorphically nondegenerate,
generic submanifold of $\bC^N$ which is minimal at some
point.  
For all $p \in M$,  there
is a unique topology on the group
$G_p$ with respect to which it is a Lie group whose Lie
algebra is
$\hol(M,p)$.
\endproclaim

This theorem follows from a slight modification of the proof
of Theorem 3.1 of Kobayashi [Ko, p. 13], by making use of
Theorem 2. Indeed, if dim$_\bR \hol(M,p) $ is finite, the
exponential of $\hol_0(M,p)$ (those vector fields in 
$\hol(M,p)$ that vanish at $p$)
generates a connected Lie group $G_p^0$ which is a
normal subgroup of
$G_p$.  One may then impose the unique topology on $G_p$
for which $G_p/G_p^0 $ is discrete. (For other applications of
this approach to groups of automorphisms see e.g.
Burns-Shnider [BS].)  However, there is a natural topology
for the group
$G_p$ obtained by regarding $G_p$ as a subspace of the space
of holomorphic mappings. Also, using Corollary of \S 1, one
can embed $G_p$ as a subgroup of the group of invertible
$(d+1)l(M)$-jets from which it inherits a topology.  In
general one does not know if these three topologies coincide.
This question will be addressed in future work [BER2].

\enddocument
\end